\documentclass[10pt,a4paper]{article}
\usepackage{amsmath, amssymb}
\usepackage{latexsym, color}
\usepackage{epsfig}

\numberwithin{equation}{section}

\newtheorem{thm}{Theorem}[section]

\newtheorem{lem}[thm]{Lemma}
\newtheorem{prop}[thm]{Proposition}

\newcommand{\qed}{\hfill \ensuremath{\square}}

\newcommand{\pf}{\noindent {\sl Proof}. \ }
\newcommand{\p}{\partial}

\newcommand{\abs}[1]{\lvert#1\rvert}

\newcommand{\eqnref}[1]{(\ref {#1})}

\newcommand{\Cbb}{\mathbb{C}}

\newcommand{\Rbb}{\mathbb{R}}

\newcommand{\Ical}{\mathcal{I}}
\newcommand{\Jcal}{\mathcal{J}}

\newcommand{\Kcal}{\mathcal{K}}

\newcommand{\Pcal}{\mathcal{P}}

\newcommand{\Scal}{\mathcal{S}}

\def\Bq{{\bf q}}

\def\BG{{\bf G}}

\def\BM{{\bf M}}

\def\BQ{{\bf Q}}

\def\BU{{\bf U}}

\def\BY{{\bf Y}}

\newcommand{\Ga}{\alpha}
\newcommand{\Gb}{\beta}

\newcommand{\Gg}{\gamma}

\newcommand{\Gk}{\kappa}

\newcommand{\Gl}{\lambda}

\newcommand{\GG}{\Gamma}

\newcommand{\beq}{\begin{equation}}
\newcommand{\eeq}{\end{equation}}

\def\Dt{{\widetilde{D}}}

\begin{document}
\title{Invariance Properties of Generalized Polarization Tensors and Design of Shape Descriptors in Three Dimensions\thanks{\footnotesize This work was
supported by the ERC Advanced Grant Project MULTIMOD--267184 and
NRF grants No. 2010-0017532.}}

\author{Habib Ammari\thanks{\footnotesize Department of Mathematics and Applications, Ecole Normale Sup\'erieure,
45 Rue d'Ulm, 75005 Paris, France (habib.ammari@ens.fr, han.wang@ens.fr).} \and
Daewon Chung\thanks{Department of Mathematics, Inha University, Incheon
402-751, Korea (chdaewon@inha.ac.kr, hbkang@inha.ac.kr).}  \and
Hyeonbae Kang\footnotemark[3]  \and Han Wang\footnotemark[2]}

\date{}
\maketitle

\begin{abstract}
We derive transformation formulas for the generalized polarization
tensors under rigid motions and scaling in three dimensions, and
use them to construct an infinite number of invariants under those
transformations. These invariants can be used as shape descriptors
for dictionary matching.
\end{abstract}



%%%%%%%%%%%%%%%%%%%%%%%%%%%%%
\section{Introduction}
%%%%%%%%%%%%%%%%%%%%%%%%%%%%%

Shape of a domain can be represented in terms of various physical
and geometric quantities such as eigenvalues, capacity and
moments. The generalized polarization tensor (GPT) is one of them.
GPTs are an (infinite) sequence of tensors associated with
inclusions (domains) and they appear naturally in the far field
expansion of the perturbation of the electrical field in the
presence of the inclusion. They are geometric quantities in the
sense that the full set of GPTs completely determines the shape of
the inclusion as proved in \cite{AK03}. This suggests that GPTs
has richer information on the shape than eigenvalues since the
full set of eigenvalues does not determine the shape uniquely
\cite{GWW}. Moreover, recent studies \cite{AGKLY, AKLZ} show that
we can use a first few terms of GPTs to recover a good
approximation of actual shape of the inclusion. Even the topology
(the number of components) can be recovered. GPTs have been used
for imaging diametrically small inclusions and computation of
effective properties of dilute composites. We refer to
\cite{AK:2007} and references therein for these applications. It
is worth mentioning that the notion of GPTs has been used not only
for imaging but also for invisibility cloaking \cite{AKLL}.

When the domain is transformed by a rigid motion or a dilation,
the corresponding GPTs change according certain rules and it is
possible to construct as combinations of GPTs invariants under
these transformations. This property makes GPTs suitable for the
dictionary matching problem. The dictionary matching problem is to
identify the object in the dictionary when the target object is
identical to one of the objects in the dictionary up to shifting,
rotation, and scaling. The standard method of dictionary matching
is to construct invariants, called shape descriptors, under rigid
motions and scaling, and to compare those invariants, and a common
way to construct such invariants uses the moments \cite{Hu}. In
the recent paper \cite{ABGJKW:pre} new invariants are constructed
using GPTs in two dimensions: it is shown that a first few terms
of GPTs and hence invariants of the target can be computed using
the measurement of the multi-static response matrix, and then the
dictionary matching technique is applied for target
identification. Viability of the method is demonstrated by
numerical experiments.

It is the purpose of this paper to extend results of
\cite{ABGJKW:pre} to construct invariants using GPTs in three
dimensions. In fact, using contracted GPTs (CGPT), which are
harmonic combinations of GPTs, we are able to construct an
infinite number of shape descriptors which are invariant under
rigid motions and scaling. Since the method of above mentioned
paper uses the complex structure of the two dimensional space, it
can not be applied to three dimensions. Using transformation
formulas of spherical harmonics under rigid motions, which can be
found in \cite{SR:1973} for example, we are able to derive
transformation formulas obeyed by CGPTs. We then use these
formulas and method of registration (see \cite{ZF}) to construct
invariants which can be used for target shape description and
position and orientation tracking \cite{ABGJKW:pre, tracking}.

%%%%%%%%%%%%%%%%%%%%%%%%%%%%%%%%%%%%%%%%%%%%%%%%%%%%%%%%%%%%%%
\section{Neumann-Poincar\'{e} operator and CGPT}\label{sec:Inv}
%%%%%%%%%%%%%%%%%%%%%%%%%%%%%%%%%%%%%%%%%%%%%%%%%%%%%%%%%%%%%%%

For a given bounded domain $D$ in $\Rbb^3$ with the Lipschitz boundary,
the Neumann-Poincar\'{e} operator $\Kcal_D$ for a density function
$\phi\in L^2(\partial D)$ and its $L^2$-adjoint $\Kcal_D^{\ast}$ are
defined in the principal value by
\begin{align*}
  \Kcal_D[\phi](x)&=\frac{1}{4\pi}\int_{\partial D}\frac{\langle
    y-x,\nu(y)\rangle}{|x-y|^3}\phi(y)
  d\sigma (y),\quad x\in\partial D\,,\\
  \Kcal^{\ast}_D[\phi](x)&=\frac{1}{4\pi}\int_{\partial D}
  \frac{\langle x-y,\nu(x)\rangle}{|x-y|^3}\phi(y)d\sigma(y),\quad
  x\in\partial D\,.
\end{align*}
Here $\langle \cdot,\cdot\rangle$ denotes the scalar product and
$\nu(x)$ the unit outward normal vector along the boundary at $x$.
Let $\Gl$ be a real number such that $|\Gl| > 1/2$. The
generalized polarization tensor (GPT) $M_{\Ga\Gb}$ for
multi-indices $\Ga=(\Ga_1, \Ga_2, \Ga_3)$ and $\Gb=(\Gb_1, \Gb_2,
\Gb_3)$ associated with $\Gl$ and $D$ is defined by \beq
M_{\Ga\Gb}(\Gl, D)= \int_{\p D} y^\Gb (\lambda
I-\Kcal^{\ast}_D)^{-1} [\nu \cdot \nabla y^\Ga] d\sigma. \eeq Here
$y^\Ga=y_1^{\Ga_1} y_2^{\Ga_2} y_3^{\Ga_3}$. Throughout this paper
$\Gl$ is fixed, so we use $M_{\Ga\Gb}(D)$ for $M_{\Ga\Gb}(\Gl,
D)$.

The contracted GPTs (CGPT) are harmonic combinations of GPTs. To
be more precise, let $Y_n^m$, $-n\le m \le n$, be the (complex)
spherical harmonic of homogeneous degree $n$ and order $m$, {\it
i.e.},
$$
Y^m_n(\theta,\varphi)=(-1)^m\bigg[\frac{2n+1}{4\pi}\frac{(n-m)!}{(n+m)!}\bigg]^{1/2}
e^{im\varphi}\Pcal^m_n(\cos\theta)\,,\quad-n\leq m\leq n\,,
$$
where $\Pcal_n^m$ are the associated Legendre polynomials of degree $n$ and order $m$.
If
$$
r^n Y_n^m(\theta, \varphi) = \sum_{|\alpha|=n} a_{\alpha}^{mn} x^\alpha,
$$
then CGPT $M_{mnkl}$ is defined by
\beq\label{CGPT2}
 M_{nmlk}= \sum_{|\alpha|=n,|\beta|=l}  \overline{a_{\alpha}^{mn}} a_{\beta}^{kl} M_{\alpha\beta}, \quad m,n,k,l=1,2,\cdots.
\eeq
In other words, we have
\beq\label{CGPT}
M_{nmlk} =\int_{\p D}r^l_y\overline{Y^k_l(\theta_y,\varphi_y)}(\lambda I-\Kcal^{\ast}_D)^{-1}\bigg[\frac{\p}{\p\nu}r^n_yY^m_n(\theta_y,\varphi_y)\bigg|_{\p D}\bigg](y)d\sigma(y),
\eeq
where $y=r_y(\cos\varphi_y\sin\theta_y,\sin\varphi_y\sin\theta_y,\cos\theta_y)$.

We now show that the operator $\Kcal_D^*$ is invariant under the
shift, scaling, and rotation of the domain $D$, which will be a
crucial fact in study of invariance properties of CGPTs in
following sections.  Let $z$ be a point in $\Rbb^3$, $s$ a
positive number, and $R$ a $3\times 3$ orthogonal matrix, and
define
\begin{itemize}
\item Shift: $T_zD=D_z=\{x_z=x+z \,|\, x\in D\}$\,; \item
  Scaling: $sD=\{sx\,|\, x\in D\}$ \,;\item Rotation:
  $R(D)=D_R=\{x_R=Rx\,|\,x\in D\}$\,.
\end{itemize}
For a function $\phi$ defined on $\p D$, define
$$
  \phi^z(y_z)=\phi(y), \quad \phi^s(sy)=\phi(y), \quad
  \phi^R(y_R) =\phi(y)\,.
$$
Then one can see easily that the following invariance holds:
\begin{align}
  \Kcal^{\ast}_{Dz}[\phi^z](x_z)&= \Kcal^{\ast}_{D}[\phi](x)\, \label{shiftform} \\
  \Kcal^{\ast}_{sD}[\phi^s](sx)&=\Kcal^{\ast}_D[\phi](x)\, \label{scaleform} \\
  \Kcal^{\ast}_{D_R}[\phi^R](x_R)&=\Kcal^{\ast}_{D}[\phi](x). \label{rotform}
\end{align}
In fact, for example, we have by the simple change of variables $y=\widetilde{y}/s$
\begin{align*}
  \Kcal^{\ast}_{sD}[\phi^s](sx)&=\frac{1}{4\pi}\int_{\partial (sD)}\frac{\langle sx-\widetilde{y},v(sx)\rangle}{|sx-\widetilde{y}|^3}\phi^s({y})d\sigma(\widetilde{y})\\
  &=\frac{1}{4\pi}\int_{\partial D}\frac{s\langle x-y,v(sx)\rangle}{s^3|x-y|^3}\phi(y)s^2 d\sigma(y)=\Kcal^{\ast}_D[\phi](x)\,.
\end{align*}
The other two relation can be seen similarly.

%%%%%%%%%%%%%%%%%%%%%%%%%%%%%%%%%%
\section{The MSR matrix and the CGPT block matrix}\label{sec:MSR}
%%%%%%%%%%%%%%%%%%%%%%%%%%%%%%%%%%

We will investigate invariance of the block matrix consisting of
CGPTs rather than individual GPTs $M_{\alpha\beta}$. The way to
construct the block matrix can be seen most clearly using the
multi-static response (MSR) matrix.

Let $\{x_r\}_{r=1}^{N}$ and $\{x_s\}_{s=1}^N$ be a set of electric potential point detectors and electric point sources. Let $u_s(x)$ be the solution to the transmission problem in the presence of inclusion $D$
\begin{equation}
  \left\{\begin{array}{ll}
      \nabla\cdot(1+(\Gk-1)\chi_D)\nabla u_s(x)=\delta_{x_s}(x)\,, & x\in\Rbb^2\setminus\p D, \\
      u_s(x)|_+=u_s(x)|_-\,, & x\in \p D, \\
      \nu(x)\cdot(\nabla u_s)|_+=\Gk\nu(x)\cdot(\nabla u_s)|_-\,, & x\in \p D, \\
      u_s(x)-\GG_s(x)=\mathcal{O}(|x|^{-2}) & |x-x_s| \rightarrow\infty\,,
    \end{array}\right.\label{trans_pro}
\end{equation}
where the notation $\phi|_{\pm}(x)$ means the limit $\lim_{t\rightarrow 0}\phi(x\pm t\nu(x))$, $\Gk= (2\Gl+1)/(2\Gl-1)$, and
$$
\GG_s(x) = \GG(x-x_s)= -\frac{1}{4\pi} \frac{1}{|x-x_s|}.
$$
Without the inclusion $D$ the solution is $\GG_s(x)$.
The MSR matrix $\mathbf{V}$ is the matrix of differences of electric potentials with and without the conductivity inclusions, and its $rs$-components are defined by
$$V_{rs}=u_s(x_r)-\GG_s(x_r),\quad 1\leq r,s\leq N\, .$$

Let $S_D$ be the single layer potential associated with $D$:
$$\mathcal{S}_D[\phi](x):=\int_{\partial D}\GG(x-y)\phi(y)d\sigma(y)\,,~x\in\Rbb^3\,.$$
It is known that the solution $u_s$ to \eqnref{trans_pro} can be represented as
\beq
u_s(x)=\GG_s(x)+\mathcal{S}_D[\phi_s](x)
\eeq
where $\phi_s\in L^2(\p D)$ solves
$$(\lambda I-\Kcal^{\ast}_D)[\phi_s]=\frac{\p\GG_s}{\p\nu}\bigg|_{\p D}\,.$$
So, the MSR matrix is given by
$$
V_{sr} = \int_{\p D}\GG(x_r-y)(\lambda I-\Kcal^{\ast}_D)^{-1}\bigg[\frac{\p\GG_s}{\p \nu}\bigg|_{\p D}\bigg](y)d\sigma(y)\,.
$$

Let $x=r_x(\cos\varphi_x\sin\theta_x,\sin\varphi_x\sin\theta_x,\cos\theta_x)$ and $y=r_y(\cos\varphi_y\sin\theta_y,\sin\varphi_y\sin\theta_y,\cos\theta_y)$ in spherical coordinates, and suppose that $r_y<r_x$. Then it is well known (see \cite{SR:1973, N:2001, AH:2012} for example) that
\beq
  \GG(x-y) =-\sum_{l=0}^{\infty}\sum_{k=-l}^l\frac{1}{2l+1}Y^k_l(\theta_x,\varphi_x)\overline{Y^k_l(\theta_y,
    \varphi_y)}\frac{r_y^l}{r_x^{l+1}} \, .\label{Gamma_sp}
\eeq
So, by assuming the inclusion $D$ is away from the sources, we have from \eqnref{CGPT}
\begin{align}
  V_{rs} &=\int_{\p D} \bigg(\sum_{l=0}^{\infty} \sum_{k=-l}^{l} \frac{1}{2l+1}Y^k_l(\theta_{x_r},\varphi_{x_r}) \overline{Y^k_l(\theta_y,\varphi_y)} \frac{r_y^l}{r_{x_r}^{l+1}} \bigg) \nonumber \\
  &\quad\times (\lambda I-\Kcal^{\ast}_D)^{-1} \bigg[ \frac{\p}{\p\nu} \sum_{n=0}^{\infty} \sum_{m=-n}^{n} \frac{1}{2n+1} Y^m_n(\theta_y,\varphi_y) \overline{Y^m_n(\theta_{x_s},\varphi_{x_s})} \frac{r_y^{n}}{r_{x_s}^{n+1}} \bigg|_{\p D} \bigg](y) d\sigma(y) \nonumber \\
  &=\sum_{l=1}^\infty \sum_{k=-1}^l \sum_{n=1}^\infty \sum_{m=-n}^n \frac{1}{(2l+1)(2n+1)} \frac{1}{r^{l+1}_{x_r}} Y^k_l(\theta_{x_r},\varphi_{x_r}) \nonumber \\
  &\quad\times \bigg(\int_{\p D}r^l_y\overline{Y^k_l(\theta_y,\varphi_y)}
  (\lambda I-\Kcal^{\ast}_D)^{-1}\bigg[\frac{\p}{\p\nu}r^n_yY^m_n(\theta_y,\varphi_y)\bigg|_{\p D}\bigg](y)d\sigma(y)\bigg)\frac{1}{r^{n+1}_{x_s}}\overline{Y^m_n(\theta_{x_s},\varphi_{x_s})} \nonumber \\
  &=\sum_{l=1}^\infty \sum_{k=-l}^l \sum_{n=1}^\infty \sum_{m=-n}^n\frac{1}{(2l+1)(2n+1)} \frac{1}{r^{l+1}_{x_r}}Y^k_l(\theta_{x_r},
  \varphi_{x_r})
  M_{nmlk} \frac{1}{r^{n+1}_{x_s}}\overline{Y^m_n(\theta_{x_s},\varphi_{x_s})}. \label{msrone}
\end{align}
Here we have used the fact that $M_{nmlk}=0$ for $n=0$ or $l=0$.

We now introduce matrix $\BM_{ln}$ by
\beq\label{mlnkm}
(\mathbf{M}_{ln})_{km}:=M_{nmlk}, \quad -l\leq k\leq l, \ -n\leq m\leq n.
\eeq
We emphasize that the dimension of $\BM_{ln}$ is $(2l+1) \times (2n+1)$. We also define $1 \times (2l+1)$ and $1 \times (2n+1)$ matrices (vectors) $\BY_{rl}$ and $\BY_{sn}$ by
\begin{align}
  (\BY_{rl})_{k}& :=\frac{Y^k_l(\theta_{x_r},\varphi_{x_r})}{(2l+1)r^{l+1}_{x_r}}, \quad -l \le k \le l \,, \\
  \quad(\BY_{sn})_{m} & = \frac{Y_m^n(\theta_{x_s},\varphi_{x_s})}{(2n+1)r^{n+1}_{x_s}}, \quad -n \le m \le n\,.
\end{align}
Then \eqnref{msrone} yields, after truncating terms corresponding to $l>K$ and $n>K$ for some integer $K$,
\beq
V_{sr}=\sum_{l,n=1}^K\BY_{rl}\mathbf{M}_{ln} \BY_{sn}^*,
\eeq
where $*$ denotes the Hermitian transpose ($A^*=(\overline{A})^t$, $t$ for transpose).

Let $N$ be the number of (coincident) source transmitters and receivers and $K$ be the truncation order as before. We further define
\begin{equation}
  \BM:= \begin{bmatrix}
      \BM_{11} & \BM_{12} & \cdots &\BM_{1K}\\
      \BM_{21} & \BM_{22} & \cdots &\BM_{2K}\\
      \vdots          & \vdots          &\ddots  &\vdots\\
      \BM_{K1} & \BM_{K2} & \cdots &\BM_{KK}\\
    \end{bmatrix}, \quad \BY:= \begin{bmatrix}
      \BY_{11} & \BY_{12} &\cdots & \BY_{1K}\\
      \BY_{21} & \BY_{22} &\cdots & \BY_{2K}\\
      \vdots          & \vdots          &\ddots & \vdots         \\
      \BY_{N1} & \BY_{N2} &\cdots & \BY_{NK} \\
    \end{bmatrix} .
\end{equation}
Then $\BM$ and $\BY$ are $(K^2+2K) \times (K^2+2K)$ and $N \times (K^2+2K)$ matrices,
respectively, and the MSR matrix $\mathbf{V}$ can be written as
\beq
\mathbf{V}= \BY \BM \BY^*
\eeq
after truncation.  We call $\BM_{ln}$ and $\BM$ a CGPT matrix and the CGPT block matrix of order $K$, respectively.

\begin{prop}\label{M:Hermitian}
  The CGPTs matrix $\BM$ is Hermitian, i.e.,
  $\BM=\overline{\BM^t}$. Furthermore, the matrices $\BM_{nn}$
  are invertible for all $n\geq 1$.
\end{prop}
\pf
Define the coefficients $a^{kl}_{\alpha}$ and $a^{mn}_{\alpha}$ so that
  $$
  r^kY^{k}_{l}(\theta,\varphi)=\sum_{|\beta|=l} a^{kl}_{\beta} y^{\beta}
  \quad \mbox{and} \quad r^nY^{m}_n(\theta,\varphi)=\sum_{|\alpha|=n}a^{mn}_{\alpha}y^{\alpha}
  $$
Then from \eqnref{CGPT2} and the symmetry property of $M_{\alpha,\beta}$ on the coefficients of harmonic
  polynomials \cite[Theorem 4.10]{AK:2007}, we have
$$
M_{nmlk}=\sum_{|\alpha|=n,|\beta|=l} \overline{a^{kl}_{\beta}}a^{mn}_{\alpha}M_{\alpha\beta}=\overline{\sum_{|\alpha|=n,|\beta|=l} a^{kl}_{\beta}\overline{a^{mn}_{\alpha}}M_{\beta\alpha}} =\overline{M_{lknm}}\,,
$$
which yields $\BM_{ln} = \BM_{nl}^*$. So $\BM$ is Hermitian.

To prove the invertibility of $\BM_{nn}$,
  it suffices to show $v^*\BM_{nn}v\neq 0$ for any $v\in\Cbb^{2n+1}$, $v \neq 0$. By definition of
  the CGPT matrix, we have
  \begin{align*}
    v^*\BM_{nn} v &= \sum_{k=1}^{2n+1}\sum_{m=1}^{2n+1}v_k^*M_{nmnk}v_m\\
    &=\int_{\p D}\overline{\sum_{k=1}^{2n+1}v_k r^n_y
      Y^k_n(\theta_y,\varphi_y)}(\lambda
    I-\Kcal^{\ast}_D)^{-1}\bigg[\frac{\p}{\p\nu}\sum_{m=1}^{2n+1}v_m
    r^n_yY^m_n(\theta_y,\varphi_y)\bigg|_{\p D}\bigg](y)d\sigma(y)\,.
  \end{align*}
Note that $\sum_{k=1}^{2n+1}v_k r^n_y Y^k_n(\theta_y,\varphi_y)$ is a harmonic
  polynomial, and we introduce the coefficients $\gamma_\alpha$ so that
  \begin{align*}
    \sum_{\abs\alpha=n}\gamma_\alpha y^\alpha = \sum_{k=1}^{2n+1}v_k r^n_y Y^k_n(\theta_y,\varphi_y).
  \end{align*}
And write $\gamma_\alpha = \gamma_\alpha^1 + i\gamma_\alpha^2$, with
  $\gamma_\alpha^1$ and $\gamma_\alpha^1$ the real and imaginary part,
  respectively. Remark that both $\gamma_\alpha^1$ and $\gamma_\alpha^2$ are
  coefficients of some harmonic polynomials. Then using the fact that $v^*\BM_{nn}v$
  is a real number, we get
  \begin{align*}
    v^*\BM_{nn}v = \sum_{\substack{\abs\alpha=n\\ \abs\beta=n}} \gamma_\alpha^*
    M_{\alpha\beta}\gamma_\beta =\sum_{\substack{\abs\alpha=n\\ \abs\beta=n}}
    \gamma_\alpha^1 M_{\alpha\beta}\gamma_\beta^1 + \gamma_\alpha^2
    M_{\alpha\beta}\gamma_\beta^2,
  \end{align*}
  while this quantity is strictly positive if $\lambda>1/2$, and strictly negative if
  $\lambda\leq 1/2$ by the positivity of $M_{\alpha\beta}$ \cite[Theorem 4.11]{AK:2007}. This completes the proof. \qed

%%%%%%%%%%%%%%%%%%%%%%%%%%%%%%%%%%%%%%%%%%%%%%%%%%%%%%%%%%%%%%
\section{Transformation formulas for the CGPT matrix}\label{sec:trans}
%%%%%%%%%%%%%%%%%%%%%%%%%%%%%%%%%%%%%%%%%%%%%%%%%%%%%%%%%%%%%%

In this section, we derive transformation formulas of the CGPT
matrix $\BM_{ln}$ defined in \eqnref{mlnkm} and the CGPT block
matrix under rigid motions and dilation. These formulas play an
essential role in finding the invariant under these
transformations in the sequel.

\subsection{Scaling}
\label{sec:CGPT_scaling}

We consider first the scaling of $\BM_{ln}$. Let us denote
$$\phi_{D,nm}(y)=(\lambda-\Kcal^{\ast}_{D})^{-1}\left[\frac{\p}{\p \nu}r^n_{y}Y^n_m(\theta_{y},\varphi_{y})\bigg|_{\p D}\right](y)\,.$$
Using the change of variables $y_s=sy$ ($r_{y_s}=sr_y,~\theta_{y_s}=\theta_y,~\varphi_{y_s}=\varphi_{y}$), we obtain
\begin{align*}
  M_{nmlk}(D)&=\int_{\p D}r^l_yY^l_k(\theta_y,\varphi_y)\phi_{D,nm}(y)d\sigma (y)\\
  &=\int_{\p (sD)}s^{-l} r^l_{y_s}Y^l_k(\theta_{y_s},\varphi_{y_s})\phi_{D,nm}(\frac{1}{s} y_s)s^{-2}d\sigma(y_s)\\
  &=s^{-(l+2)}\int_{\p (sD)} r^l_{y_s}Y^l_k(\theta_{y_s},\varphi_{y_s})\phi^s_{D,nm}(y_s)d\sigma(y_s)\,.
\end{align*}
Since
$$
\left\langle \nu(y),\nabla(r^n_yY^n_m(\theta_y,\varphi_y))\right\rangle=s^{-n+1}\left\langle \nu(y_s),\nabla(r^n_{y_s}Y^n_m(\theta_{y_s},\varphi_{y_s})\right\rangle,
$$
it follows from \eqnref{scaleform} that
\begin{align*}
(\lambda I-\Kcal^{\ast}_{sD})[\phi^s_{D,nm}](y_s) &=(\lambda I-\Kcal^{\ast}_D) [\phi_{D,nm}](y) \\
& = \left\langle \nu(y),\nabla(r^n_{y}Y^n_m(\theta_{y},\varphi_{y})\right\rangle\\
&=s^{-n+1}\left\langle \nu(y_s),\nabla(r^n_{y_s}Y^n_m(\theta_{y_s},\varphi_{y_s})\right\rangle\,,
\end{align*}
and hence
$$
\phi^s_{D,nm}(y_s)=s^{-n+1}\phi_{sD,nm}(y_s)\,.
$$
Thus, we have
$$
M_{nmlk}(D)=s^{-(l+n+1)}\int_{\p sD}r^l_{y_s}Y^l_k(\theta_{y_s},\varphi_{y_s})\phi_{sD,nm}(y_s)d\sigma(y_s)=s^{-(l+n+1)}M_{nmlk}(sD)\,.
$$

\begin{lem}[Scaling]
  For any positive integers $l,n$ and scaling parameter $s>0$, the following holds:
\begin{equation}
    \BM_{ln}(sD)=s^{l+n+1}\BM_{ln}(D)\,.\label{scalingM}
\end{equation}
\end{lem}

\subsection{Shifting}
\label{sec:CGPT_shifting}

To deal with a shifting of $\BM_{ln}$, we need the translation of the regular
spherical harmonics $r^n Y_n^m(\theta,\varphi)$. For $y=(r,\theta,\varphi)$,
$z=(r_z,\theta_z,\varphi_z)$, and $y_z=y+z=(r',\theta',\varphi')\,,$ we have the
following expression of the translation of the regular spherical harmonic,
\beq\label{shiftRSH}
  r'^nY_n^m(\theta',\varphi')
  =\sum_{(\nu,\mu)}^{(n,m)} C_{\nu\mu nm} r_z^{n-\nu}Y_{n-\nu}^{m-\mu}(\theta_z,\varphi_z)r^{\nu}Y_{\nu}^{\mu}(\theta,\varphi)\,,
\eeq
where
$$
C_{\nu\mu nm}=\bigg[\frac{4\pi(2n+1)(n-m)!(n+m)!}{(2n-2\nu+1)(2\nu+1)(n-\nu-m+\mu)!(n-\nu+m-\mu)!(\nu-\mu)!(\nu+\mu)!}\bigg]^{1/2}\,.
$$
Here we use special summation notation:
$$
\sum_{(\nu,\mu)}^{(n,m)} = \sum_{\nu=0}^n\sum_{\mu=\max(-\nu,\nu-n+m)}^{\mu=\min(\nu,-\nu+n+m)}.
$$
We refer to \cite{SR:1973} for a proof of formula
\eqnref{shiftRSH}. Using \eqnref{shiftform} and (\ref{shiftRSH}),
we obtain
\begin{align}
  M_{nmlk}(D_z)&=\int_{\p D} r'^l\overline{Y^k_l(\theta',\varphi')}(\Gl I -\Kcal^{\ast}_D)^{-1}\bigg[\frac{\p}{\p \nu}r'^nY^m_n(\theta',\varphi')\bigg|_{\p D}\bigg](y)d\sigma(y)\nonumber\\
  &=\int_{\p D}\sum_{(i,j)}^{(l,k)} C_{ijlk} r_z^{l-i}\overline{Y_{l-i}^{k-j}(\theta_z,\varphi_z)r^{i}Y_{i}^{j}(\theta,\varphi)} \nonumber\\
  &\quad \times (\Gl I -\Kcal^{\ast}_D)^{-1} \bigg[\frac{\p}{\p \nu}\sum_{(\nu,\mu)}^{(n,m)} C_{\nu\mu nm} r_z^{n-\nu}Y_{n-\nu}^{m-\mu}(\theta_z,\varphi_z)r^{\nu}Y_{\nu}^{\mu}(\theta,\varphi)\bigg|_{\p D}\bigg](y)d\sigma(y)\nonumber\\
  &=\sum_{(i,j)}^{(l,k)}\sum_{(\nu,\mu)}^{(n,m)} C_{ijlk} r_z^{l-i}\overline{Y^{k-j}_{l-i}(\theta_z,\varphi_z)}
  M_{ij\nu\mu}(D) C_{\nu\mu nm} r^{n-\nu}_zY^{m-\mu}_{n-\nu}(\theta_z,\varphi_z)\label{F:shift}\,.
\end{align}

For a pair of integers $i$ and $l$, define $(2l+1) \times (2i+1)$ matrix
$\BG_{li}= \BG_{li}(z)$ by
\beq\label{compG}
(\BG_{li})_{kj}=G_{ijlk} :=\left\{\begin{array}{ll}
        C_{ijlk} r_z^{l-i}Y^{k-j}_{l-i}(\theta_z,\varphi_z)\quad & \textrm{if } \max(-i,i-l+k)\leq j\leq\min(i,-i+l+k),\\
        0\quad & \textrm{otherwise}
      \end{array}\right.
\eeq
for $- l \leq k \leq l$. Then, from the equation (\ref{F:shift}) and taking the fact that $\BM_{ln}=0$ if $n=0$ or $l=0$ into account, we immediately have the following property.

\begin{lem}[Shifting] For any positive integer $l,n$ and the shifting parameter $z$, the following holds:
  \begin{equation}\label{shiftM}
    \BM_{ln}(D_z)=\sum_{i=1}^{l}\sum_{\nu=1}^{n}\overline{\BG_{li}(z)}
    \BM_{i\nu}(D) \BG_{n\nu}(z)^t\,.
  \end{equation}
\end{lem}

To gain better understanding of the matrix $\BG_{li}$, we compute $\BG_{21}$. Indeed it will play a role in Section \ref{sec:Inva}. $\BG_{21}$ is given by
$$
\BG_{21}=\begin{bmatrix}
    G_{1,-1,2,-2} & G_{1,0,2,-2} & G_{1,1,2,-2} \\
    G_{1,-1,2,-1} & G_{1,0,2,-1} & G_{1,1,2,-1} \\
    G_{1,-1,2,0} & G_{1,0,2,0} & G_{1,1,2,0} \\
    G_{1,-1,2,1} & G_{1,0,2,1} & G_{1,1,2,1} \\
    G_{1,-1,2,2} & G_{1,0,2,2} & G_{1,1,2,2} \\
  \end{bmatrix}
\,.
$$
For $k=-2,-1,0,1,2$ ($l=2$ and $i=1$), the conditions $\max(-i,i-l+k)\leq j\leq \min(i,-i+l+k)$ can be written as follows:
\begin{align*}
  k=-2: ~\max(-1,-1+k)\leq j\leq \min(1,1+k)&~\Rightarrow~j=-1\\
  k=-1:~\max(-1,-1+k)\leq j\leq \min(1,1+k)&\Rightarrow j=-1,0\\
  k=0:~\max(-1,-1+k)\leq j\leq \min(1,1+k)&\Rightarrow j=-1,0,1\\
  k=1:~\max(-1,-1+k)\leq j\leq \min(1,1+k)&\Rightarrow j=0,1\\
  k=2:~\max(-1,-1+k)\leq j\leq \min(1,1+k)&\Rightarrow j=1\,,
\end{align*}
Then using (\ref{compG}) and the definition of the spherical harmonics, we have
\begin{align}
  \BG_{21}&=\begin{bmatrix}
      G_{1,-1,2,-2} & 0            & 0 \\
      G_{1,-1,2,-1} & G_{1,0,2,-1} & 0 \\
      G_{1,-1,2,0}  & G_{1,0,2,0}  & G_{1,1,2,0} \\
      0             & G_{1,0,2,1}  & G_{1,1,2,1} \\
      0             & 0            & G_{1,1,2,2} \\
    \end{bmatrix} \nonumber \\
    &= \begin{bmatrix}
      C_{1,-1,2,-2} r_z Y^{-1}_1(\theta_z,\varphi_z) & 0                                           & 0 \\
      C_{1,-1,2,-1} r_z Y^{0}_1(\theta_z,\varphi_z)  & C_{1,0,2,-1} r_z Y^{-1}_1(\theta_z,\varphi_z) & 0 \\
      C_{1,-1,2,0} r_zY^{1}_1(\theta_z,\varphi_z)   & C_{1,0,2,0} r_zY^{0}_1(\theta_z,\varphi_z)   & C_{1,1,2,0} r_z Y^{-1}_1(\theta_z,\varphi_z)\\
      0  & C_{1,0,2,1} r_zY^{1}_1(\theta_z,\varphi_z)   & C_{1,1,2,1} r_zY^{0}_1(\theta_z,\varphi_z) \\
      0  & 0  & C_{1,1,2,2} r_z Y^{1}_1(\theta_z,\varphi_z) \\
    \end{bmatrix} \nonumber \\
  &= \begin{bmatrix}
      -\sqrt{5}r_z\sin\theta_z e^{-i\varphi_z}           & 0                                                & 0 \\
      \sqrt{5}r_z\cos\theta_z                           & -\sqrt{\frac{5}{2}}r_z\sin\theta_z e^{-i\varphi_z} & 0 \\
      \sqrt{\frac{5}{6}}r_z\sin\theta_z e^{i\varphi_z} & \sqrt{\frac{20}{3}}r_z\cos\theta_z                                 & -\sqrt{\frac{5}{6}}r_z\sin\theta_z e^{-i\varphi_z} \\
      0    & \sqrt{\frac{5}{2}}r_z\sin\theta_z e^{i\varphi_z} & \sqrt{5}r_z\cos\theta_z \\
      0  & 0    & \sqrt{5}r_z\sin\theta_z e^{i\varphi_z}\\
    \end{bmatrix} \,. \label{G21}
\end{align}
On the other hand, from the equation (\ref{compG}), one can easily see that the $G_{njnk}$ is non-zero only when $k=j$, and
$$
G_{nknk}=C_{nknk} Y_0^0(\theta_z,\varphi_z)=\sqrt{4\pi}Y_0^0(\theta_z,\varphi_z)=1\,.
$$
Thus, $\BG_{nn}$ is the identity $(2n+1) \times (2n+1)$ matrix.
So, \eqnref{shiftM} yields, for instance,
\beq
  \BM_{21}(D_z)=\overline{\BG_{21}(z)} \BM_{11}(D) + \BM_{21}(D) \,.
\eeq

\subsection{Rotation}
\label{sec:CGPT_rotation}

Rotations in three dimensions may be described in many different ways. Among them we use the Euler angles
which can be conveniently used to represent the rotation formula for spherical harmonics. The rotation $R$ is given by
$$
R= \begin{bmatrix}
    \cos\Gg & -\sin\Gg & 0 \\
    \sin\Gg & \cos\Gg & 0 \\
    0 & 0 & 1 \\
  \end{bmatrix}
  \begin{bmatrix}
    \cos\Gb & 0 & -\sin\Gb \\
    0 & 1 & 0 \\
    \sin\Gb & 0 & \cos\Gb \\
    \end{bmatrix}
  \begin{bmatrix}
    \cos\Ga & -\sin\Ga & 0 \\
    \sin\Ga & \cos\Ga & 0 \\
    0 & 0 & 1 \\
   \end{bmatrix}\,.
$$
That is, we rotate by an angle $\Ga$ about $z$-axis, and by an angle
$\Gb$ about the new $y$ axis, and finally by an angle $\Gg$ about the new
$z$-axis. Let $R$ be a rotation matrix. Since the homogeneous polynomials and Laplace operator
are invariant under rotation, $Y^m_n(R\xi)$ is also spherical harmonic of degree $n$, moreover
$Y^m_n(R\xi)$ can be written as follows:
\beq\label{rotspherical}
Y^m_n(R\xi)=\sum_{m'=-n}^n\rho_{n}^{m',m}Y^{m'}_n(\xi)\,,
\eeq
where
\begin{equation} \label{Dwigner}
  \rho_n^{m',m}=e^{im'\Gg}d_n^{m',m}(\Gb)e^{im\Ga}\,,
\end{equation}
Here,
\begin{align*}
  d_n^{m',m}(\Gb)
  &=[(n+m')!(n-m')!(n+m)!(n-m)!]^{1/2} \\
  & \quad \times \sum_s \frac{(-1)^{m'+m+s}(\cos\frac{\Gb}{2})^{2(n-s)+m-m'}(\sin\frac{\Gb}{2})^{2s-m+m'}}{(n+m-s)!s!(m'-m+s)!(n-m'-s)!}\,,
\end{align*}
where the sum is over values $s$ such that the factorials are
nonnegative:
$$\max(0,m-m')\leq s\leq \min(n-m',n+m)\,.$$
A proof of \eqnref{rotspherical} can be found in \cite{SR:1973}.

We want to find the transformation formula $M_{nmlk}$, {\it i.e.},
$M_{nmlk}(D_R)$. We have
\begin{align}
  M_{nmlk}(D_R)&=\int_{\p D_R}r^l\overline{Y^k_l(R\xi)}(\Gl I-\Kcal^{\ast}_{D_R})^{-1}\bigg[\frac{\p}{\p\nu}r^nY^m_n(R\xi)\bigg|_{\p D_R}\bigg](Ry)d\sigma(Ry) \nonumber \\
  &=\int_{\p D}r^l\sum_{k'=-l}^l\overline{\rho_{l}^{k',k}Y^{k'}_l(\xi)}(\Gl I-\Kcal^{\ast}_{D})^{-1}\bigg[\frac{\p}{\p\nu}r^n\sum_{m'=-n}^n
  \rho_{n}^{m',m}Y^{m'}_n(\xi)\bigg|_{\p D}\bigg](y)d\sigma(y) \nonumber \\
  &=\overline{\mathbf{q}^k_l}
  \BM_{ln}(\mathbf{q}^m_n)^t\,, \label{rotCGPT}
\end{align}
where
$$
\Bq_n^m = (\rho_n^{-n,m}, \ldots, \rho_n^{n,m}).
$$

Let for each positive integer $n$
  \begin{equation}
  \mathbf{Q}_n= \BQ_n(R):=\begin{bmatrix}
        \rho^{-n,-n}_{n} & \rho^{-n+1,-n}_{n} & \cdots & \rho^{n,-n}_{n} \\
        \rho^{-n, -n+1}_{n} & \rho^{-n+1, -n+1}_{n} & \cdots & \rho^{n, -n+1}_{n} \\
        \cdots & \cdots & \ddots & \cdots \\
        \rho^{-n,n}_n & \rho^{-n+1, n}_n & \cdots & \rho^{n,n}_n \\
      \end{bmatrix} \,.\label{def:Q}
  \end{equation}
The rotation matrix $\BQ_n$ is called Wigner $D$-matrix and known to be is unitary (see \cite{SR:1973}).

\begin{lem}[Rotation] For a unitary matrix $R$ the following relation holds:
  \begin{equation}
    \BM_{ln}(D_R)=\overline{\BQ_l(R)} \BM_{ln}(D)\BQ(R)^t \, . \label{rotationM}
  \end{equation}
\end{lem}

\subsection{CGPT block matrices and its properties}\label{sec:CGPT}

In addition to the CGPT block matrix, we define two block matrices: Let $K$ be the truncation order, $s$ be a scaling parameter, $z$ shifting factor, and $R$ rotation, and define
\beq
\BG(R):= \begin{bmatrix}
      \BG_{11} & 0 & \cdots & 0 \\
      \BG_{21} & \BG_{22} & \cdots &0 \\
      \vdots          & \vdots          &\ddots  &\vdots\\
      \BG_{K1} & \BG_{K2} & \cdots &\BG_{KK}\\
    \end{bmatrix}, \quad
\BQ(s, R):= \begin{bmatrix}
      s\BQ_{1} & 0 & \cdots & 0\\
      0 & s^2 \BQ_{2} & \cdots &0\\
      \vdots          & \vdots          &\ddots  &\vdots\\
      0 & 0 & \cdots & s^K \BQ_{K}\\
    \end{bmatrix},
\eeq
where $\BG_{ln}=\BG_{ln}(z)$ and $\mathbf{Q}_n=\BQ_n(R)$ are defined by (\ref{compG}) and (\ref{def:Q}), respectively. Then we have the following theorem
\begin{thm}
Let $T_z$ be the shift by $z$, $T^s$ scaling by $s$, and $R$ a unitary matrix. Then the CGPT block matrix $\BM$ satisfies
\begin{equation}\label{finalform}
  \BM(T_zT^sR (D))=s \overline{\BG(z)} \overline{\BQ(s,R)} \BM(D) \mathbf{Q}(s,R)^t\BG(z)^t \,.
\end{equation}
\end{thm}
\pf
We have from \eqnref{shiftM} that
$$
\BM(T_zT^sR (D))= \overline{\BG}\BM(T^sR (D)) \BG^t .
$$
We also have from \eqnref{scalingM} and \eqnref{rotationM} that
\begin{align*}
 \BM(T^sR (D)) &= s \begin{bmatrix}
      s^2 \BM_{11}(R(D)) & \cdots & s^{1+K} \BM_{1K}(R(D)) \\
            \vdots & \ddots  &\vdots\\
      s^{K+1} \BM_{K1}(R(D)) & \cdots & s^{2K} \BM_{KK} (R(D))\\
    \end{bmatrix} \\
    & = s \begin{bmatrix}
      s^2 \overline{\BQ_1} \BM_{11}(D) \BQ_1^t & \cdots & s^{1+K} \overline{\BQ_1} \BM_{1K}(D) \BQ_K^t \\
      \vdots          &\ddots  &\vdots\\
      s^{K+1} \overline{\BQ_K} \BM_{K1}(D) \BQ_1^t & \cdots & s^{2K} \overline{\BQ_K} \BM_{KK} (D) \BQ_K^t \\
    \end{bmatrix}.
\end{align*}
So we obtain \eqnref{finalform}. \qed

%%%%%%%%%%%%%%%%%%%%%%%%%%%%%%%%%%%%%%%%%%%%%%%%%%%%%%%%%%%%%%
\section{Transform invariant shape descriptors}\label{sec:Inva}
%%%%%%%%%%%%%%%%%%%%%%%%%%%%%%%%%%%%%%%%%%%%%%%%%%%%%%%%%%%%%%%

In this section, we construct the invariants using CGPT matrices
under shifting, scaling, and rotation.

Let $B$ be a reference domain and $D$ be the one obtained by rotating $B$ by $R$, scaling by $s$, and shifting by $z$, {\it i.e.},
\beq
D=T_zT^sR(B).
\eeq
Since $\BQ_n$ is unitary, we have $\overline{\mathbf{Q}_{n}}^{-1} = \mathbf{Q}_{n}^t$. One can easily see that $\BM_{11}$ and $\overline{\mathbf{Q}_{1}}$ commute, {\it i.e.},
\beq
\BM_{11} \BQ_1 = \BQ_1 \BM_{11}.
\eeq
Since $\BG_{nn}=\mathbf{I}_{2n+1}$ ($(2n+1) \times (2n+1)$ identity
matrix), we have
\begin{align*}
  \BM_{11}(D)&=s^3\overline{\mathbf{Q}_1}\BM_{11}(B)\mathbf{Q}_{1}^t = s^3 \BM_{11}(B) \,,\\
  \BM_{21}(D)&=s^3\overline{\BG_{21}\mathbf{Q}_1}\BM_{11}(B)\mathbf{Q}^t_1+s^4\overline{\mathbf{Q}_2}
  \BM_{21}(B)\mathbf{Q}_1^t=
  \overline{\BG}_{21}\BM_{11}(D)+s^4 \overline{\mathbf{Q}_2}\BM_{21}(B)\mathbf{Q}_1^t\,.
\end{align*}
Define
$$
\BU_D=\BM_{21}(D)\BM_{11}(D)^{-1} \quad\mbox{and}\quad \BU_B=\BM_{21}(B)\BM_{11}(B)^{-1}.
$$
Then we have
\beq\label{eq:UD_UB}
  \BU_D = \overline{\BG_{21}} + s\overline{\BQ_2}\BU_B\BQ_1^t.
\eeq

\subsection{Invariance by registration}
\label{sec:recalage}

We first recall that the matrices $\BG_{ln}$ and $\BQ_n$ are determined by the shift factor $z$ and rotation $R$, respectively. Moreover,
in view of (\ref{G21}), we may write $\overline{\BG}_{21}$ in terms of rectangular coordinates as
\begin{equation}
    \overline{\BG_{21}(z)} =\sqrt{5} \begin{bmatrix}
      -(z_1+z_2i)          & 0       & 0 \\
      z_3                &-\sqrt{\frac{1}{2}}(z_1+z_2i) & 0 \\
      \sqrt{\frac{1}{6}}(z_1-z_2i)& \sqrt{\frac{4}{3}}z_3                       & -\sqrt{\frac{1}{6}}(z_1+z_2i) \\
      0                          & \sqrt{\frac{1}{2}}(z_1-z_2i)& z_3 \\
      0                          & 0                          & z_1-z_2i \\
    \end{bmatrix} \,.
  \label{eq:Gb_21}\end{equation}

We present here a method of registration to construct invariants. This method is based on a linear mapping $u: M_{5 \times 3}(\Cbb) \rightarrow\Cbb^3$ ($M_{5 \times 3}(\Cbb)$ is the collection of $5 \times 3$ complex matrices) which
satisfies
\beq\label{umapcond2}
u(\overline{\BG_{21}(z)})=z
\eeq
and
\begin{align}
  \label{eq:u_map_cond}
  u(\overline{\BQ_2(R)}\BU\BQ_1(R)^t) = R u(\BU) \ \mbox{ for any $\BU \in M_{5 \times 3}(\Cbb)$ and any rotation } R.
\end{align}
Such a linear mapping does exist and an example is
\begin{align}
  \label{eq:u_example}
  u(\BU) = \frac{1}{\sqrt 5}
  \begin{bmatrix}
    -\frac{3}{10} \BU_{11} + \frac{\sqrt 3}{10\sqrt 2} \BU_{31} - \frac{3}{10\sqrt 2}
    \BU_{22} + \frac{3}{10\sqrt 2} \BU_{42} - \frac{\sqrt 3}{10\sqrt 2} \BU_{33} +
    \frac{3}{10} \BU_{53}\\
    i(\frac{3}{10} \BU_{11} + \frac{\sqrt 3}{10\sqrt 2} \BU_{31} + \frac{3}{10\sqrt 2}
    \BU_{22} + \frac{3}{10\sqrt 2} \BU_{42} + \frac{\sqrt 3}{10\sqrt 2} \BU_{33} +
    \frac{3}{10} \BU_{53})\\
    \frac{3}{10} \BU_{21} + \frac{\sqrt 3}{5} \BU_{32} + \frac{3}{10} \BU_{43}
  \end{bmatrix}.
\end{align}
One can easily check that this linear transformation satisfies \eqnref{umapcond2} using \eqnref{eq:Gb_21}.
We check that it satisfies \eqnref{eq:u_map_cond} through symbolic computations using Matlab.
Since the computation is lengthy, we omit the detail here.

In the following we use the shorthand $u_D=u(\BU_D)$. Applying $u$ to \eqnref{eq:UD_UB}, we have
\beq\label{eq:uD_uB}
  u_D = z + sR u_B,
\eeq so $u_D$ can be used as a registration point (see \cite{ZF}).
We emphasize that since $\BU_D=\BM_{21}(D)\BM_{11}(D)^{-1}$, $u_D$
can be computed using CGPTs of $D$.

The first invariant we introduce is
\beq
\mathbf{\mathcal{J}}_{ln}(D):=\BM_{ln}(T_{-u_D} D)= \overline{\BG_l(-u_D)} \BM_{ln}(D) \BG_l(-u_D)^t \, \quad l, n = 1,2,\ldots.
\eeq

\begin{prop}
  For all indices $l,n$ the quantity $\Jcal_{ln}$ is invariant by translation:
  \begin{align}
    \label{eq:translation_invar_J}
    \Jcal_{ln}(T_zD) = \Jcal_{ln}(D) \ \mbox{ for any } z\in\Rbb^3 .
  \end{align}
\end{prop}
\pf Let $D_z=T_z D$. By definition of $\Jcal_{ln}$, and using
formula \eqref{shiftM}, we have
  \begin{align*}
    \Jcal_{ln}(T_z D) &= \BM_{ln}(T_{-u_{D_z}}D_z) = \BM_{ln}(T_{z-u_{D_z}}D)\\
    &= \overline{\BG_l(z-u_{D_z})}\BM_{ln}(D) \BG_n(z-u_{D_z})^t.
  \end{align*}
The relation \eqref{eq:uD_uB} says that $u_{D_z} = z+u_D$, and hence we have
  $$ \Jcal_{ln}(T_z D) = \overline{\BG_l(-u_D)}\BM(D) \BG_n (-u_D)^t = \Jcal_{ln}(D).$$
This completes the proof. \qed

In the following we denote $\Dt=D-u_D$, then
\begin{lem}
Let $\Dt=T_{-u_D} D$.  For any indices $l,n$ and scaling parameter $s>0$ and rotation $R$, we have
\beq\label{eq:fJ_scl_rot}
    \Jcal_{ln}(sR(D)) = \Jcal_{ln}(sR(\Dt)) = \BM_{ln}(sR(\Dt)) = s^{l+n+1}\BM_{ln}(R(\Dt)).
\eeq
In particular, we have
\beq\label{eq:fJ_scl_rot2}
    \Jcal_{ln}(sD)= s^{l+n+1} \Jcal_{ln}(D),
\eeq
and
\beq\label{eq:fJ_scl_rot3}
    \Jcal_{ln}(R(D))= \overline{\BQ_l(R)} \Jcal_{ln}(D) \BQ_n(R)^t.
\eeq
\end{lem}
\pf
Since $\Jcal_{ln}$ is translation invariant, we have for any $z\in\Rbb^3$:
  \begin{align*}
    \Jcal_{ln}(sR(D)) = \Jcal_{ln}(T_{sRz}sR(D)) = \Jcal_{ln}(sR(T_zD)).
  \end{align*}
Then by taking $z=-u_D$, we obtain the first identity in \eqnref{eq:fJ_scl_rot}. Moreover, we have from \eqref{eq:uD_uB}
$$
u_{\Dt} = -u_D + u_D=0.
$$
So we have $u(sR\Dt) = sR u_{\Dt} = 0$ , which implies the second identity in \eqnref{eq:fJ_scl_rot}. The third one is \eqnref{scalingM}. By taking $R=I$, we have \eqnref{eq:fJ_scl_rot2}. \eqnref{eq:fJ_scl_rot3} follows from the definition of $\Jcal_{ln}$ and \eqnref{rotationM}.
\qed

In particular, it can be seen from \eqref{eq:fJ_scl_rot} that all $\Jcal_{nn}$ are
invertible. So we define the second invariant:
\begin{align}
  \label{eq:invar_S}
  \Scal_{ln}(D) := \Jcal_{nn}(D)^{-1} \Jcal_{nl}(D)\Jcal_{ll}(D)^{-1} \Jcal_{ln}(D).
\end{align}
It is worth emphasizing that $\Scal_{ln}(D)$ is a square matrix of dimension $(2n+1)$.

\begin{prop}
  For any indices $l,n$, the quantity $\Scal_{ln}(D)$ is translation and scaling
  invariant:
  \begin{align}
    \label{eq:fS_trl_scl_invar}
    \Scal_{ln}(T_zsD) = \Scal_{ln}(D) \ \mbox{ for any } z\in\Rbb^3, \mbox{ and } s>0.
  \end{align}
  Moreover, for any rotation $R$:
  \begin{align}
    \label{eq:fS_rot}
    \Scal_{ln}(R(D)) = \overline{\BQ_{n}(R)} \Scal_{ln}(D) \BQ_{n}(R)^t.
  \end{align}
\end{prop}
\pf Since $\Jcal_{ln}$'s are translation invariant, so are
$\Scal_{ln}$'s. Scaling invariance of $\Scal_{ln}$ follows from
\eqnref{eq:fJ_scl_rot2}. Formula \eqnref{eq:fS_rot} can be seen
using \eqnref{eq:fJ_scl_rot3} as follows:
\begin{align*}
\Scal_{ln}(R(D)) &= \Jcal_{nn}(R(D))^{-1} \Jcal_{nl}(R(D))\Jcal_{ll}(R(D))^{-1} \Jcal_{ln}(R(D)) \\
&= \overline{\BQ_n} \Jcal_{nn}(D)^{-1}\BQ_n^t \overline{\BQ_n} \Jcal_{nl}(D) \BQ_l^t \overline{\BQ_l} \Jcal_{ll}(D)^{-1} \BQ_l^t \overline{\BQ_l} \Jcal_{ln}(D) \BQ_n^t \\
&= \overline{\BQ_n} \Jcal_{nn}(D)^{-1} \Jcal_{nl}(D) \Jcal_{ll}(D)^{-1} \Jcal_{ln}(D) \BQ_n^t \\
&= \overline{\BQ_n} \Scal_{ln}(D) \BQ_n^t .
\end{align*}
This completes the proof. \qed

Since the Frobenius norm of a matrix remains unchanged after the multiplication by a
unitary matrix, so finally we define the shift descriptor
\begin{align}
  \label{eq:fI_tsr_invar}
  \Ical_{ln}(D) = \| \Scal_{ln}(D)\|_F
\end{align}
which is clearly invariant by any rotation, scaling and translation.

% Let $\sqrt{\frac{1}{6}}(\mathbf{U})_{11}=u_1+iu_2\q\textrm{and}\q\sqrt{\frac{1}{3}}(\mathbf{U})_{12}=u_3\,.$
% Then $u=z+c$, for some constant independent of the unknown shifting factor $z\,.$
% We first define the translation invariant quantity as follows. $$\mathbf{\mathcal{J}}_{ln}(D)=\BM_{ln}(D_{-u})\,.$$
% We define the scaling and rotation invariant quantity as follows
% $$\mathbf{\mathcal{S}}_{ln}(D)=\mathbf{\mathcal{J}}_{nn}^{-1}(D)\mathbf{\mathcal{J}}^t_{ln}(D)\mathbf{\mathcal{J}}_{ll}^{-1}(D)\mathbf{\mathcal{J}}_{ln}(D)\,.$$
% and
% $$\mathbf{\mathcal{I}}_{ln}(D)=\|\mathbf{\mathcal{S}}_{ln}(D)\|_{F}\,,$$
% where $\|\cdot\|_{F}$ denotes the Frobenius matrix norm, which is
% invariant under rotations. As we constructed, the quantity
% $\mathcal{I}_{ln}$ are invariant under shifting, scaling, and
% rotation for all indices $l$ and $n$.

\section{Concluding remarks}
In this paper we  have constructed new shape descriptors in three
dimensions which are invariant under translation, rotation, and
scaling. These shape descriptors can be used to efficiently
identify a target using a dictionary of precomputed CGPTs data.
They can be also used for tracking the position and orientation of
a mobile three-dimensional target from multistatic measurements.


\begin{thebibliography}{99}

\bibitem{ABGJKW:pre}
  H. Ammari, T. Boulier, J. Garnier, W. Jing, H. Kang and H. Wang, Target identification
    using dictionary matching of generalized polarization tensors,
    arXiv:1204.3035.

\bibitem{tracking} H. Ammari, T. Boulier, J. Garnier, H. Kang and H. Wang,
Tracking of a mobile target using generalized polarization
tensors, submitted.

\bibitem{ADKL:pre}
  H. Ammari, Y. Deng, H. Kang, and H. Lee, Reconstruction of inhomogeneous
  conductivities via generalized polarization tensors,
  arXiv:1211.4495.

\bibitem{AGKLY} H. Ammari, J. Garnier, H. Kang, M. Lim, and S. Yu,
Generalized polarization tensors for shape description, Numerische
Math., to appear.

\bibitem{AK03} H. Ammari and H. Kang, Properties of the generalized polarization tensors,
Multiscale Modeling and Simul., 1 (2003), 335--348.

\bibitem{AK:2007}
  H. Ammari and H. Kang, \emph{Polarization and Moment Tensors: with Applicatios to Inverse Problems and Effective
  Medium Theory}, Vol. 162, Springer-Verlag, 2007.

\bibitem{AKLL} H. Ammari, H. Kang, H. Lee, and M. Lim, Enhancement of near cloaking using generalized polarization
 tensors Vanishing Structures. Part I: The conductivity problem, Comm. Math. Phys., to appear.

\bibitem{AKLZ} H. Ammari, H. Kang, M. Lim, and H. Zribi, The generalized polarization
tensors for resolved imaging. Part I: Shape reconstruction of a
conductivity inclusion, Math. Comp., 81 (2012), 367--386.

\bibitem{AH:2012}
  K. Atkinson and W. Han, \emph{Spherical Harmonics and Approximations on the Unit Sphere: An Introduction},
  Lecture Notes in Mathematics, Vol. 2044, Springer-Verlag, New York, 2012.

\bibitem{GWW} C. Gordon, D.L. Webb, and S. Wolpert, One cannot hear the shape of a drum,
Bull. Amer. Math. Soc., 27 (1992), 134--138.

\bibitem{Hu} M.-K. Hu, Visual pattern recognition by moment invariants, IRE Trans. Inform. Theory, 8 (1962), 179--187.

\bibitem{N:2001}
  J.-N N\'{e}d\'{e}lec, \emph{Acoustic and Eletromagnetic Equations. Integral Representations for Harmonic
  Problems}, Applied Mathematical Science, Vol. 144, Springer-Verlag, New York, 2001.

\bibitem{SR:1973}
  E.O. Steinborn and K. Ruedenberg, Rotation and translation of regular and irregular solid spherical harmonics,
  Adv. Quantum Chem., 7 (1973), 1--81.

\bibitem{ZF} B. Zitova and J. Flusser, Image registration methods:
a survey, Imag. Vision Comput., 21 (2003), 977--1000.
\end{thebibliography}
\end{document}